\documentclass[12pt]{article}
\usepackage{amssymb,amsmath,graphics}
\newtheorem{thm}{Theorem}
\newtheorem{lemma}[thm]{Lemma}
\newtheorem{prop}[thm]{Proposition}
\newtheorem{cor}[thm]{Corollary}
\newenvironment{pf}
{\begin{trivlist} \item \noindent{\sc Proof. }} {{\hfill $\Box$}
\end{trivlist}}
\newenvironment{pfof}[1]
{\begin{trivlist} \item \noindent{\sc Proof of\ #1. }} {{\hfill
$\Box$} \end{trivlist}}

\newenvironment{fakethm}
{\begin{trivlist} \item \em}{\end{trivlist}}

\newcommand{\zt}{{\mathbb Z}^2}

\newcommand{\prob}{{\mathbb P}}

\newcommand{\lng}{\mbox{\rm long}}

\newcommand{\cl}[1]{\langle #1 \rangle}

\newcounter{mycount}
\newenvironment{mylist}{\begin{list}{{\rm (\roman{mycount})}}%
{\usecounter{mycount}\itemsep 0pt}}{\end{list}}

\newcommand{\dof}{\bf}

\title{Slow Convergence in Bootstrap Percolation}
\author{
Janko Gravner\thanks{Funded in part by NSF Grant DMS-0204376
and the Republic of Slovenia's Ministry of Science Program P1-285}
 \and
Alexander E. Holroyd\thanks{Funded in part by an NSERC (Canada)
Discovery Grant, and by Microsoft Research} }
\date{May 1, 2007}

\begin{document}
\maketitle

\begin{abstract}In the bootstrap percolation model, sites in an $L$ by $L$ square
are initially infected independently with probability $p$. At
subsequent steps, a healthy site becomes infected if it has at least
2 infected neighbours.  As $(L,p)\to(\infty,0)$, the probability
that the entire square is eventually infected is known to undergo a
phase transition in the parameter $p\log L$, occurring
asymptotically at $\lambda=\pi^2/18$ \cite{h-boot}. We prove that
the discrepancy between the critical parameter and its limit
$\lambda$ is at least $\Omega((\log L)^{-1/2})$.  In contrast, the
critical window has width only $\Theta((\log L)^{-1})$.  For the
so-called modified model, we prove rigorous explicit bounds which
imply for example that the relative discrepancy is at least $1\%$
even when $L=10^{3000}$.  Our results shed some light on the
observed differences between simulations and rigorous asymptotics.
\renewcommand{\thefootnote}{}
\footnote{\hspace{-2em}{\bf Key words:} bootstrap percolation,
cellular automaton, metastability, finite-size scaling, crossover}
\footnote{\hspace{-2em}{\bf 2000 Mathematics Subject
Classifications:} Primary 60K35; Secondary 82B43}
\end{abstract}

\section{Introduction}

The {\dof standard bootstrap percolation model} on the square lattice
$\zt$ is defined as follows.  For any set $K\subseteq\zt$ we define
$${\cal B}(K):=K\cup\big\{x\in\zt: \#\{y\in K: \|x-y\|_1=1\}\geq
2\big\},$$ and
$$\cl{K} :=\lim_{t\to \infty} {\cal B}^t(K),$$ where ${\cal B}^t$ denotes
the $t$-th iterate of the function ${\cal B}$.  The set $\cl{K}$ is
the final set of infected sites if we start with $K$ infected.

Now fix $p\in(0,1)$ and let $W$ be a random subset of $\zt$ in
which each site is included independently with probability $p$;
more formally let $\prob=\prob_p$ be the product measure with
parameter $p$ on $\Omega=\{0,1\}^{\zt}$, and define the random
variable $W=W(\omega):=\{x\in \zt:\omega(x)=1\}$ for
$\omega\in\Omega$.  We say that a set $K\subset\zt$ is {\dof
internally spanned} if $\cl{K\cap W} = K$.  For $L\geq 1$ denote
the square $R(L):=\{1,\ldots ,L\}^2\subset\zt$.  The main object
of interest is the function
$$I(L)=I(L,p):=\prob_p\big( R(L) \text{ is internally spanned}\big).$$

A central result is the following from \cite{h-boot}.
\begin{fakethm}
{\bf Theorem (phase transition, \cite{h-boot})} Consider the
standard bootstrap percolation model. As $L\to\infty$ and $p\to 0$
simultaneously we have
\begin{equation}
\begin{array}{rl}
\text{if}\quad \liminf \,p \log L > \lambda \quad&\text{then}\quad I(L,p)\to 1;\\
\text{if}\quad \limsup \,p \log L < \lambda \quad&\text{then}\quad
I(L,p)\to 0,
\end{array}\label{h-phase}
\end{equation}
where $\lambda:={\pi^2}/{18}$.
\end{fakethm}

Surprisingly, predictions for the asymptotic threshold $\lambda$
based on simulation differ greatly from the rigorous result.  For
example, in \cite{a-s-a} the estimate $\lambda=0.245\pm 0.015$ is
reported (based on simulation of squares up to size $L=28800$),
whereas in fact $\lambda=\pi^2/18=0.548311\cdots$.  This apparent
discrepancy between theory and experiment has been investigated
using partly non-rigorous methods in \cite{dlbd,dld,stauffer}.  Our
aim is to provide some rigorous understanding of the phenomenon. Our
main result is the following strengthening of the first assertion in
\eqref{h-phase}.

\begin{thm}[slow convergence]
\label{slow} Consider the standard bootstrap percolation model.
There exists $c>0$ such that, if $L\to\infty$ and $p\to 0$
simultaneously in such a way that
$$p\log L > \lambda - \frac{c} {\surd \log L}\, ,$$
where $\lambda=\pi^2/18$, then
$$I(L,p)\to 1.$$
\end{thm}
(The condition in Theorem \ref{slow} may be equivalently expressed
as $p\log L>\lambda - c'\surd p$, for a different constant $c'$).
Thus, the convergence of the critical value of the parameter $p \log
L$ to its limit $\lambda$ is very slow, with an asymptotic
discrepancy of at least $c/\surd (\log L)$. (In order to halve the
latter quantity, $L$ must be raised to the 4th power).

On the other hand, the window over which $I$ changes from near $0$
to near $1$ is much smaller -- roughly  $\mbox{constant}/\log L$.
The precise statement depends on whether we vary $p$ or $L$, as
follows.

For fixed $L$, and $\alpha\in(0,1)$, define
$p_\alpha=p_\alpha(L):=\sup\{p: I(L,p)\le \alpha\}$. Since $I(L,p)$
is continuous and strictly increasing in $p$, we have that
$p_\alpha$ is the unique value such that $I(L,p_\alpha)=\alpha$. The
following was proved in \cite{balogh-bollobas-sharp} using a general
result from \cite{friedgut-kalai}.

\begin{fakethm}\sloppy
{\bf Theorem ($p$-window, \cite{balogh-bollobas-sharp})} Consider
the standard bootstrap percolation model.  For any fixed
$\epsilon\in(0,1)$, we have
\begin{equation}\label{p-win}
p_{1-\epsilon}\log L-p_\epsilon\log L=O\bigg(\frac{\log\log
L}{\log L}\bigg)=O\big(p_{1/2}\log p_{1/2}^{-1}\big) \quad\text{as
} L\to\infty.
\end{equation}
\end{fakethm}

More precise estimates on the size of the window are available if we
instead vary $L$.  An upper bound was proved in
\cite{aizenman-lebowitz}. Here we use similar methods to obtain
matching upper and lower bounds. Since $I(L,p)$ is not necessarily
monotone in $L$, we define for fixed $p$ and $\alpha\in(0,1)$:
$\underline{L}_\alpha=\underline{L}_\alpha(p):= \min\{L: I(L,p)\ge
\alpha\}$ and $\overline{L}_\alpha=\overline{L}_\alpha(p):= \max\{L:
I(L,p)\le \alpha\}$. Thus the interval
$[\underline{L}_\epsilon,\overline{L}_{1-\epsilon}]$ contains all
those $L$ for which $I(L,p)\in[\epsilon,1-\epsilon]$.

\begin{thm}[$L$-window]\sloppy
\label{l-win} Consider the standard bootstrap percolation model.
For any fixed $\epsilon\in(0,1/5)$, we have
$$p\log \overline{L}_{1-\epsilon}-p\log \underline{L}_\epsilon =
\Theta(p)=\Theta\Big(1/{\log \overline{L}_{1/2}}\Big)
\quad\text{as } p\to 0.$$ Indeed, for $p$ sufficiently small we
have
$$p\log \overline{L}_{1-\epsilon}-p\log \underline{L}_\epsilon \in
[C_-p,C_+p],$$ where $C_\pm=C_\pm(\epsilon)= (1/2\pm o(1))\log
\epsilon^{-1}$ as $\epsilon\to 0$.
\end{thm}

The {\dof modified bootstrap percolation model} is a variant of
the standard model in which we replace the update rule ${\cal B}$
with
$${\cal B}_M(K):=K\cup\big\{x\in\zt: \{x+e_i,x-e_i\}\cap
K\neq\emptyset \text{ for each of }i=1,2 \big\}$$ (here $e_1:=(1,0)$
and $e_2:=(0,1)$ are the standard basis vectors), and define
$\cl{\cdot}_M$, internally spanned, and $I_M(L,p)$ accordingly.
\begin{fakethm}
{\bf Theorem (\cite{h-boot})}  For the modified bootstrap
percolation model, \eqref{h-phase} holds with threshold
$\lambda_M:=\pi^2/6$.
\end{fakethm}

\begin{thm}\label{mod}
Theorem \ref{l-win} and \eqref{p-win} hold also for the modified
model.
\end{thm}

In the case of Theorem \ref{slow} we establish the following
stronger version with an explicit error bound.
\begin{thm}[explicit bound]
\label{explicit} For the modified model, if $p\leq 1/10$ and
$$p\log L \geq \lambda_M - \sqrt{2p}+{\cal E}(p),\quad\text{then}\quad
I_M(L,p)\geq 1/2,$$ where $\lambda_M=\pi^2/6$ and ${\cal
E}(p):=1.8 p\log p^{-1} + 2 p$.
\end{thm}
One may deduce rigorous numerical bounds such as the following.
\begin{cor}
\label{very-explicit} Consider the modified model.  We have
$p_{1/2}\log L < 0.98 \,\lambda_M$ when $L=10^{500}$, and
$p_{1/2}\log L < 0.99 \,\lambda_M$ when $L=10^{3000}$.
\end{cor}
\begin{pf}
Take respectively $p=0.0014$ and $p=0.0002356$ in Theorem
\ref{explicit}.
\end{pf}

\subsubsection*{Remarks}

Aside from their mathematical interest, bootstrap percolation models
have been applied to a variety of physical problems (see e.g.\
\cite{adler-vis}), and as tools in the study of other models (e.g.\
\cite{martinelli,fontes-schonmann-sidoravicius,frobose}).

Several interesting attempts have been made to understand the
discrepancy between simulation results (e.g.\ those of \cite{a-s-a})
and the rigorous results in \cite{h-boot}; see e.g.\
\cite{adler-vis,dlbd,dld,stauffer}. The present work is believed to
be the first fully rigorous progress in this direction. In
\cite{stauffer} it is estimated that $p_{1/2}\log L$ may become
close to $\lambda=\pi^2/18$ only beyond about $L=10^{20}$ (the data
given in \cite{a-s-a} support a similar conclusion). Current
simulations extend only to about $L=10^5$.  A length scale of about
$L=10^{10}$ is relevant to some physical applications.  Thus it is
important to understand this issue in more detail.

In particular, it would be of interest to determine the asymptotic
behaviour of (say) $\lambda-p_{1/2}\log L$ as $L\to\infty$.  Theorem
\ref{slow} gives only a lower bound of $\Omega((\log L)^{-1/2})$. In
\cite{stauffer} simulation data is fitted to $p_{1/2}\log
L=\pi^2/18-0.45(\log L)^{-0.2}$.  In \cite{dld}, computer
calculations together with a heuristic argument lead to the estimate
$p_{1/2}\log L=\pi^2/6-3.67 (\log L)^{-0.333}$ for the modified
model.


The phenomenon of a critical window whose width is asymptotically
much smaller than its distance from a limiting value has been proved
in other settings including integer partitioning problems
\cite{borgs-chayes-pittel}, but contrasts with more familiar models
such as random graphs \cite{luczak}.

\subsubsection*{Outline of Proofs}

The idea behind the phase transition result \eqref{h-phase} from
\cite{h-boot} is as follows.  We expect the square $R(L)$ to be
internally spanned if and only if it contains at least one
internally spanned square of side $B \gg 1/p$, since with high
probability this will grow indefinitely in the presence of a
random background of density $p$.  Such a square is sometimes
called a nucleation centre or critical droplet. Therefore
 the critical regime should be roughly at $L^2 I(B)\approx 1$, i.e.\
 $\log L\approx (-\log I(B))/2$, and we need to estimate $I(B)$.
First consider the modified model.  One way for $R(B)$ to be
internally spanned is for every square with its bottom left corner
at $(1,1)$ to have at least one adjacent occupied site on each its
top and right faces -- then every such square will be internally
spanned (we can think of an infected square growing from $R(1)$ to
$R(B)$).  A straightforward computation shows that the probability
of this event is approximately $\exp[-2\lambda_M/p]$ where
$\lambda_M=\pi^2/6$. This argument proves the first inequality in
\eqref{h-phase} for the modified model. (The second inequality
requires a much more delicate argument - see \cite{h-boot}).

In order to prove the slow convergence result for the modified
model, Theorem \ref{explicit}, we consider other ways for a square
to be internally spanned.  One way is for every site along the main
diagonal to be occupied.  For a square of size $A<p^{1/2}$, the
latter event has higher probability than the event in the previous
paragraph, because the probability of growing by one additional row
and column is $p$ versus about $(Ap)^2$. Therefore let
$A=p^{-1/2}/2$, and suppose $R(A)$ is internally spanned by this
mechanism, while each square from $R(A)$ to $R(B)$ has occupied
sites on its faces as before.  By comparing the two growth
mechanisms, we see that, compared with the previous argument, this
increases the lower bound on $I(B)$ by a factor of least
$[p/(Ap)^2]^A=\exp[Cp^{-1/2}]$. This argument therefore proves the
analogue of Theorem \ref{slow} for the modified model. Theorem
\ref{explicit} is proved by a refinement of these ideas (see in
particular Lemmas \ref{diag} and \ref{scan}). The coefficient $\surd
2$ of $\surd p$ seems to be the best that can be achieved by this
method.

The above argument cannot work for the standard bootstrap
percolation model.  This is because an internally spanned square can
grow from a face whenever there is an occupied site within distance
2.  Thus, each additional occupied site can allow growth by two rows
or two columns, so we do not achieve sufficient saving by
considering occupied sites along the diagonal.  Instead we consider
another mechanism.  Rather than a growing square, we consider a
growing rectangle which may change shape when it encounters vacant
rows or columns. (Figure \ref{mechanism} illustrates the main idea).
We may describe such growth by means of the path traced by the
rectangle's top right corner. As noted in \cite{h-boot}, the
probability of such a growth path becomes much smaller if it
deviates far from the main diagonal (which corresponds to a growing
square). However, it turns out that if the deviations are of scale
only $p^{-1/2}$ then the ``entropy factor'' (the number of possible
deviations) outweighs the ``energy cost'' (the reduction in
probability for each path). This argument yields Theorem \ref{slow}.

\subsubsection*{Notation}
\label{sec-pre}

The following notation will be used throughout.  For integers
$a,b,c,d$ we denote the rectangle
$(a,b;c,d):=([a,c]\times[b,d])\cap\zt$, and we write for
convenience $R(m,n)=R(1,1;m,n)$ and $R(n)=R(n,n)$.  The {\dof long
side} of a rectangle is $\lng(R(a,b;c,d))=\max\{c-a+1,d-b+1\}$. A
{\dof copy} of a set $K\subseteq\zt$ is an image under an isometry
of $\zt$.  A site $x\in\zt$ is {\dof occupied} if $x\in W$.  A set
of sites is {\dof vacant} if it contains no occupied site.

It will sometimes be convenient to denote
$$q=q(p):=-\log(1-p),$$
and
$$f(z):=-\log(1-e^{-z}),$$
so that for any $K\subset \zt$,
$$\prob_p(K \text{ is not vacant})=1-(1-p)^{|K|}=\exp-f(|K|q).$$
Note that $q\geq p$, and $q\sim p$ as $p\to 0$.  The function $f$ is
positive, decreasing, and convex on $(0,\infty)$.

In Section \ref{sec-slow} we will also have occasion to consider the
functions
$$\beta(u):=\frac{u+\sqrt{u(4-3u)}}{2} \quad\text{and}\quad
g(z):=-\log \beta(1-e^{-z}).
$$
The thresholds $\lambda,\lambda_M$ arise from the integrals
\begin{equation}\label{integral}
\int_0^\infty f = \lambda_M=\frac{\pi^2}{6} \quad\text{and}\quad
\int_0^\infty g = \lambda=\frac{\pi^2}{18}
\end{equation}
 (see \cite{h-boot}).

\section{Critical Window}

In this section we present a proof of Theorem \ref{l-win}, together
with the extension to the modified model claimed in Theorem
\ref{mod}.
The following lemma from \cite{aizenman-lebowitz} is useful.
\begin{lemma}
\label{a-l} Let $R$ be a rectangle, and consider the standard or
modified model.  If $R$ is internally spanned then for every
positive integer $k\leq \lng(R)$ there exists an internally spanned
rectangle $T\subseteq R$ with $\lng(T)\in[k,2k]$.
\end{lemma}

\begin{pf}  See \cite{aizenman-lebowitz}.  \end{pf}

\begin{lemma}[comparison]
\label{comp} Consider the standard or modified model. For integers
$L\geq \ell\geq 2$ and any $p\in(0,1)$ we have
\begin{mylist}
\item
$$I(L)\geq \Big(1-e^{-I(\ell)\big(\tfrac L\ell-1\big)^2}\Big)
\left( 1-2L^2 e^{-p\ell} \right);$$
\item
$$\left( 1-2\ell^2 e^{-p(\ell/4-1)} \right) I(L) \leq I(\ell)
\bigg(\frac{2L}{\ell-1}\bigg)^2.$$
\end{mylist}
\end{lemma}

\begin{pfof}{Lemma \ref{comp}(i)}
Let $m=\lfloor L/\ell\rfloor$, and consider the $m^2$ disjoint
squares
$$S_k=R(\ell)+k\ell, \qquad k\in\{0,\ldots,m-1\}^d.$$
Let $E$ be the event that at least one of the $S_k$ is internally
spanned, and let $F$ be the event that every copy of $R(1,\ell)$ in
$R(L)$ is non-vacant.  It is straightforward to see that if $E$ and
$F$ both occur then $R(L)$ is internally spanned. Hence using the
Harris-FKG inequality (see e.g.\ \cite{g2}),
\begin{align*}
I(L)\geq \prob(E)\prob(F)\geq &
\Big( 1-(1-I(\ell))^{m^2}\Big) \Big(1-2L^2(1-p)^\ell\Big) \\
\geq & \Big(1-e^{-I(\ell)\big(\tfrac L\ell-1\big)^2}\Big) \left(
1-2L^2 e^{-p\ell} \right).
\end{align*}
\end{pfof}

\begin{pfof}{Lemma \ref{comp}(ii)}
Let $s=\lfloor \ell/2\rfloor$ and $m=\lfloor L/s\rfloor$, and
consider the $m^2$ overlapping squares
$$S_k=R(\ell)+ ks \wedge (L-\ell,L-\ell), \qquad k\in\{0,\ldots,m-1\}^2,$$
where $\wedge$ denotes coordinate-wise minimum.  Note that
$\bigcup_k S_k=R(L)$, and that the overlap between two adjacent
squares has width at least $s$. It follows that any rectangle
$T\subseteq R(L)$ with $\lng(T)\leq s$ lies entirely within one of
the $S_k$. Hence, using Lemma \ref{a-l},
\begin{align}
I(L)&\leq  \prob\left( \exists \text{ i.s. }T \subseteq R(L) \text{
with }
\lng(T)\in \big[\lfloor \tfrac s2 \rfloor,s\big] \right) \nonumber\\
&\leq  \prob \Big[\bigcup_k \big\{ \exists \text{ i.s. }T \subseteq
S_k \text{ with }
\lng(T)\in\big[\lfloor \tfrac s2 \rfloor,s\big] \big\}\Big] \nonumber\\
\label{comp-step}
 &\leq  m^2 \prob\big( \exists \text{ i.s. }T \subseteq R(\ell)
\text{ with } \lng(T)\in\big[\lfloor \tfrac s2 \rfloor,s\big] \big).
\end{align}

On the other hand, considering the event that every copy of
$R(1,\lfloor \tfrac s2 \rfloor)$ in $R(\ell)$ contains at least one
occupied site, and using the argument from the proof of part (i), we
have
$$I(\ell)\geq \prob\big( \exists \text{ i.s. }T \subseteq R(\ell) \text{ with }
\lng(T)\in\big[\lfloor \tfrac s2 \rfloor,s\big] \big) \big(1-2\ell^2
e^{-ps}\big).$$
Combining this with (\ref{comp-step}) yields the result.
\end{pfof}

\begin{pfof}{Theorem \ref{l-win}}
It follows from \eqref{h-phase} that for any $\alpha\in(0,1)$ we
have
\begin{equation}
\label{plogl}
 p\log \overline{L}_\alpha(p) \;,\; p\log \underline{L}_\alpha(p) \to \lambda
 \quad\text{as }p\to 0.
 \end{equation}
Therefore, once the first equality is proved, the second follows
immediately. To prove the first equality we will use Lemma
\ref{comp} to derive upper and lower bounds on $p\log
\overline{L}_{1-\epsilon}-p\log \underline{L}_\epsilon$.

For the upper bound, we fix $\epsilon$ and use Lemma \ref{comp}(i)
with $L=\overline{L}_{1-\epsilon}(p)$ and
$\ell=\underline{L}_\epsilon(p)$, noting that $I(L,p)\leq
1-\epsilon$ and $I(\ell,p)\geq \epsilon$.  By (\ref{plogl}), for $p$
sufficiently small (depending on $\epsilon$) we have $1-2L^2 e^{-p
\ell} \geq 1-\epsilon^2$, so we obtain for $p$ sufficiently small:
$$1-\epsilon \geq \Big(1-e^{-\epsilon
\big(\tfrac {\overline{L}_{1-\epsilon}} {\underline{L}_\epsilon}
-1\big)^2}\Big)
 (1-\epsilon^2).$$
Rearranging gives
$$\frac {\overline{L}_{1-\epsilon}} {\underline{L}_\epsilon} \leq
1+ \sqrt{\frac{1}{\epsilon} \log \frac{1+\epsilon}{\epsilon}},$$
hence
$$p\log
\overline{L}_{1-\epsilon}-p\log \underline{L}_\epsilon \leq C_+ p,$$
where
$C_+=\log\big(1+\sqrt{\epsilon^{-1}\log(\epsilon^{-1}+1)}\,\big)$
satisfies $C_+<\infty$ for all $\epsilon>0$ and
$C_+\leq(\tfrac12+o(1))\log \epsilon^{-1}$ as $\epsilon\to 0$.

For the lower bound, we fix $\epsilon$ and use Lemma \ref{comp}(ii)
with $L=\overline{L}_{1-\epsilon}(p)+1$ and
$\ell=\underline{L}_\epsilon(p)-1$, noting that $I(L,p)> 1-\epsilon$
and $I(\ell,p)< \epsilon$.  By (\ref{plogl}), we have $2\ell^2 e^{-p
(\ell/4-1)}=o(1)$ as $p\to 0$, so we obtain:
$$\big(1-o(1)\big) (1-\epsilon) \leq \epsilon
\Big(\frac{2(\overline{L}_{1-\epsilon}+1)}{\underline{L}_\epsilon-2}\Big)^2.$$
Rearranging gives
$$\frac {\overline{L}_{1-\epsilon}+1} {\underline{L}_\epsilon-2} \geq
\sqrt{\frac{(1-\epsilon)(1-o(1))}{4 \epsilon}},$$ as $p\to 0$.  For
$p$ sufficiently small we obtain
$$p\log
\overline{L}_{1-\epsilon}-p\log \underline{L}_\epsilon \geq C_- p,$$
for any $C_-(\epsilon)<\log\sqrt{(1-\epsilon)/(4\epsilon)}$.  Thus
we may take $C_->0$ for all $\epsilon<1/5$, and $C_-\geq
(\tfrac12-o(1))\log\epsilon^{-1}$ as $\epsilon\to 0$.
\end{pfof}

\section{Slow Convergence}
\label{sec-slow}

The main step in proving Theorem \ref{slow} will be the following.

\begin{prop}[nucleation centres]
\label{nuc-slow}
Consider the standard bootstrap percolation model. There exist
$p_0>0$ and $c\in(0,\infty)$ such that, for all $p<p_0$ and $B\geq
2p^{-1}$,
$$I(B,p) \geq \exp\big[-2\lambda/p+c/\surd p\big],$$
where $\lambda=\pi^2/18$.
\end{prop}

\begin{pfof}{Theorem \ref{slow}}
First suppose that $(L,p)\to(\infty,0)$ in such a way that for some
$c_1$,
\begin{equation*}
\label{l-corr} p\log L >\lambda -c_1/\surd \log L.
\end{equation*}
Then for $L$ sufficiently large we have in particular $p\log L
>\lambda/2$, hence
\begin{equation}
\label{p-corr} p\log L > \lambda -c_2\surd p,
\end{equation}
where $c_2=2c_1/\lambda$.

Therefore it is enough to prove that for some $c_2>0$, if
$(L,p)\to(\infty,0)$ satisfy \eqref{p-corr} then $I(L,p)\to 1$.
Furthermore, we may assume that we have equality in \eqref{p-corr},
since if not we may find (for $p$ sufficiently small) $p'<p$ such
that $p'\log L = \lambda -c_2\surd p'$, and then $I(L,p)\geq
I(L,p')\to 1$.  Therefore let
$$L=\exp\Big[\lambda/p - c_2/\surd p\Big] \quad\text{and}\quad
B=\lceil p^{-3} \rceil.$$ Using Lemma \ref{comp}(i),
\begin{equation}\label{meta}
I(L)\geq \Big(1-e^{-I(B)\big(\tfrac L B-1\big)^2}\Big) \left( 1-2L^2
e^{-p B} \right).
\end{equation}
 Proposition \ref{nuc-slow} and the above
definition of $L$ easily imply $L^2 e^{-p B}\to 0$ as $p\to 0$,
while
$$\log\Big[ I(B)(L/B-1)^2 \Big] \leq -2\lambda/p +c/\surd p +
2\big(\lambda/p -c_2/\surd p\big) +O(\log p^{-1}) \to 0$$ as $p\to
0$ provided $2c_2>c$.  Then \eqref{meta} gives $I(L,p)\to 1$ as
required.
\end{pfof}

\begin{figure}
\centering
\begin{picture}(160,160)
\put(27,47){$a$}\put(47,27){$a$}
\put(150,70){$b$}\put(70,150){$b$}
\resizebox{!}{2in} {\includegraphics{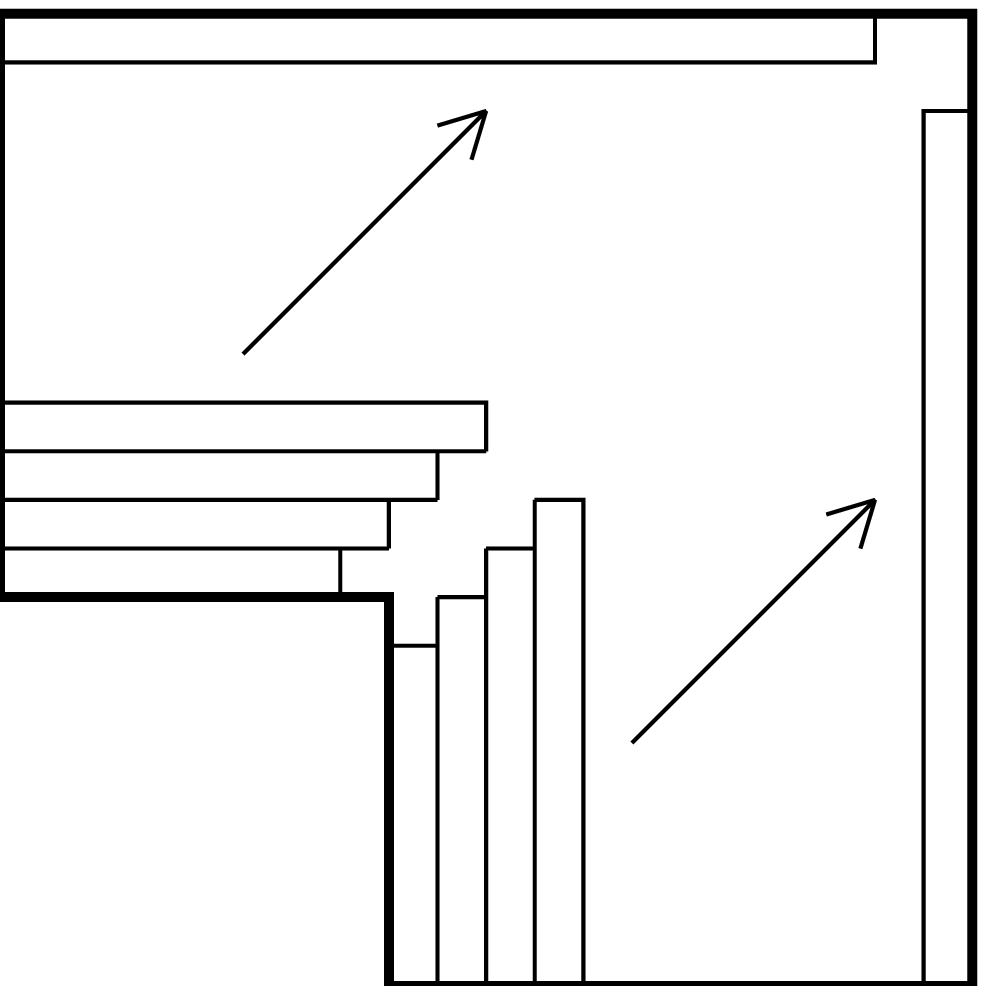}}
\end{picture}
\hfill
\begin{picture}(160,160)
\put(27,47){$a$}\put(47,27){$a$}
\put(150,70){$b$}\put(70,150){$b$}
\resizebox{!}{2in}{\includegraphics{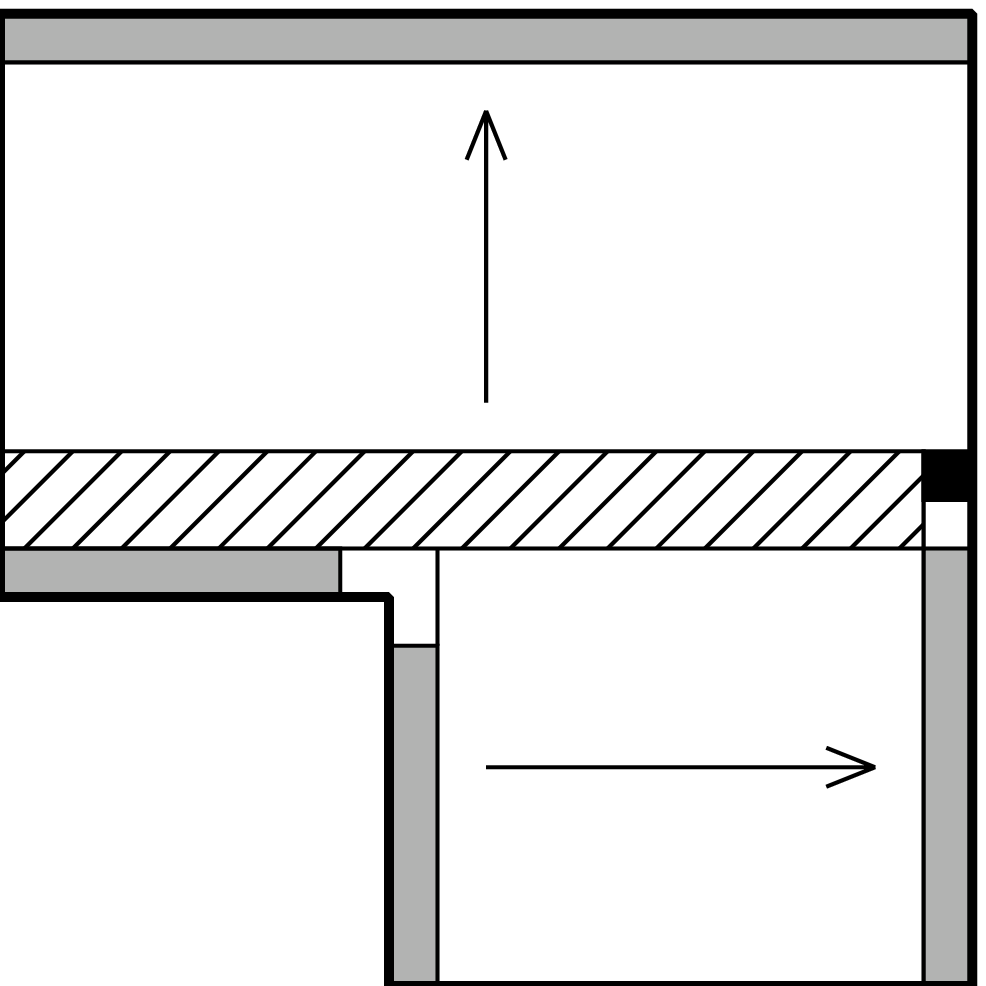}}
\end{picture}
 \caption{Two possible mechanisms for growth from $R(a)$ to $R(b)$.
 (i) The event ${\cal D}_a^b$:  no two consecutive strips are
 vacant.
 (ii) The event ${\cal J}_a^b$: the gray strips are non-vacant, the hatched region is vacant,
 the black site is occupied, and the horizontal/vertical arrows indicate no two consecutive
 vacant columns/rows respectively.}
\label{mechanism}
\end{figure}

In order to prove Proposition \ref{nuc-slow} we consider various
ways for $R(B)$ to be internally spanned.  The simplest way involves
symmetric growth starting from a corner.  We say that a sequence of
events $A_1,A_2,\ldots,A_k$ has a {\dof double gap} if there is a
consecutive pair $A_i,A_{i+1}$ neither of which occur.  For integers
$2\leq a\leq b$, let ${\cal D}_a^b$ be the event that:
\begin{align*}
\big\{ R(1,i;\;i-2,i) \text{ is not vacant}\big\}_{i=a+1,\ldots,b} &
\text{ has no double gaps, and}\\
\big\{ R(i,1;\;i,i-2) \text{ is not vacant}\big\}_{i=a+1,\ldots,b} &
\text{ has no double gaps.}
\end{align*}
See Figure \ref{mechanism}(i).  Note that if $R(a)$ is internally
spanned, and ${\cal D}_a^b$ occurs, then $R(s,t)$ is internally
spanned for some $s,t\in\{b-1,b\}$.  Indeed, it is easily seen that
we may find a sequence of internally spanned rectangles $R(i,j)$
with $|i-j|\leq 2$, starting with $R(a)$ and ending with $R(s,t)$,
with the width or the height increasing by 1 or 2 at each step.

We will also consider the following alternative growth mechanism.
For positive integers $a\leq b-4$, let ${\cal J}_a^b$ be the event
that:
\begin{align*}
R(1,a+1;\;a-1,a+1) &\text{ is not vacant},\\
R(a+1,1;\;a+1,a-1) &\text{ is not vacant},\\
\big\{ R(i,1;\;i,a+1) \text{ is not
vacant}\big\}&_{i=a+2,\ldots,b-1}
 \text{ has no double gaps}, \\
 (b,1;\;b,a+1) &\text{ is not vacant}, \\
R(1,a+2;\;b-1,a+3) &\text{ {\em is} vacant},\\
(b,a+3) &\text{ is occupied},\\
\big\{ R(1,i;\;b,i) \text{ is not vacant}\big\}&_{i=a+4,\ldots,b-1}
 \text{ has no double gaps, and} \\
R(1,b;\;b,b) &\text{ is not vacant}.
\end{align*}
See Figure \ref{mechanism}(ii).  Note again that if $R(a)$ is
internally spanned and ${\cal J}_a^b$ occurs then $R(b)$ is
internally spanned. In this case, vertical growth is stopped by the
two vacant rows, and there is a sequence of horizontally growing
internally spanned rectangles, followed by vertical growth after the
occupied site $(b,a+3)$ is encountered.

Now fix a positive integer $B$.  For positive integers
$(a_i,b_i)_{i=1,\ldots,m}$ satisfying $2\leq a_1\leq b_1\leq a_2\leq
\cdots \leq b_m\leq B$ and  $b_i-a_i\leq 4 \;\forall i$, define the
event
\begin{align*}
{\cal E}(a_1,b_1,\ldots ,a_m,b_m)&:= {\cal D}_2^{a_1} \cap
\Big(\bigcap_{i=1}^{m} {\cal J}_{a_i}^{b_i}\Big) \cap
\Big(\bigcap_{i=1}^{m-1} {\cal D}_{b_i}^{a_{i+1}}\Big) \cap {\cal
D}_{b_m}^{B-1} \\ & \cap \Big\{(1,1),(2,2),(B,1),(1,B) \text{
 are occupied}\Big\}.
\end{align*}

\begin{lemma}[properties of ${\cal E}$]{\ }\label{properties}
\begin{mylist}
\item
The various events appearing in the above definition of ${\cal
E}(a_1,\ldots ,b_m)$ are independent.
\item
If ${\cal E}(a_1,\ldots ,b_m)$ occurs then $R(B)$ is internally
spanned.
\item
For different choices of $a_1,\ldots ,b_m$, the events ${\cal
E}(a_1,\ldots ,b_m)$ are disjoint.
\end{mylist}
\end{lemma}

\begin{pf}
Property (i) is clear from the definitions of the ${\cal D}$ and
${\cal J}$ events.  Property (ii) follows from the earlier remarks
on these events: indeed the squares $R(2),R(b_1),\ldots,R(b_m),R(B)$
are all internally spanned. To see (iii), fix a configuration and
consider examining in sequence the rows $R(1,i;\;i-2,i)$ for
$i=3,4,5,\ldots$.  The presence of two consecutive vacant rows
signals an event ${\cal J}_a^b$, and determines the value of $a$,
and then if we follow the upper vacant row to the right until an
occupied site is encountered, we discover the corresponding value of
$b$.
\end{pf}

We will obtain a lower bound on the probability $R(B)$ is
internally spann\-ed by bounding the probability of each event
${\cal E}$ (for certain choices of the $a_i,b_i$), and bounding
the number of possible choices.

We start by estimating the probability of ${\cal D}_a^b$, for which
we need the following slight refinement of a result from
\cite{h-boot} (see \cite{andrews} for a much more precise result in
the same direction). Recall the function $\beta$ defined in the
introduction.
\begin{prop}[double gaps]
\label{recur} For independent events $A_1,\ldots,A_k$ whose
probabilities $u_i:=\prob(A_i)$ form an increasing or decreasing
sequence, the probability that there are no double gaps is at least
$\prod_{i=1}^k \beta(u_i)$.
\end{prop}

\begin{lemma}
\label{ineq} For $0\leq u\leq v\leq 1$ we have $u \beta(v) +(1-u)v
\geq \beta(u)\beta(v)$.
\end{lemma}

\begin{pf}
The function $h(u,v):=u \beta(v) +(1-u)v - \beta(u)\beta(v)$
satisfies $h(v,v)=0$, so it suffices to show that $h$ is decreasing
in $u$ for $u\leq v$.  But we have $\partial h/\partial u =
\beta(v)-v-\beta'(u)\beta(v)\leq 0$, by the elementary computations
$\beta'(u)\geq \beta'(v)\geq (\beta(v)-v)/\beta(v)$.
\end{pf}

\begin{pfof}{Proposition \ref{recur}}
Without loss of generality suppose the probabilities $u_i$ are
decreasing.  Let $a_k$ be the probability that the sequence
$A_1,\ldots,A_k$ has no double gaps.  Then $a_0=a_1=1$, and by
conditioning on the last two events we obtain $a_k=u_k
a_{k-1}+(1-u_k) u_{k-1} a_{k-2}$.  The result follows by induction,
using Lemma \ref{ineq} thus: $a_k \geq  [u_k \beta(u_{k-1}) +
(1-u_k) u_{k-1}] \prod_{i=1}^{k-2} \beta(u_i) \geq \prod_{i=1}^{k}
\beta(u_i)$.
\end{pfof}

Recall the function $g$ from the introduction, and write for $a\leq
b$,
$$G_a^b=G_a^b(p):=\exp \Big[\textstyle  - \sum_{i=a}^{b-1} g(iq)\Big].$$
\begin{lemma}[diagonal growth]
\label{diag-bound}
$$\prob_p({\cal D}_a^b) \geq (G_{a-1}^{b-1})^2$$
\end{lemma}

\begin{pf}
Immediate from Proposition \ref{recur} and the definitions of ${\cal
D}_a^b$ and $g$.
\end{pf}

Next we estimate the relative cost of a ${\cal J}$-event.
\begin{lemma}[deviation cost]
\label{jig-bound}
 Fix positive constants $c_-<c_+$.  For any $p\in(0,1/2)$ and
$a\leq b-4$ satisfying $a,b\in [c_-/p,c_+/p]$, we have
$$\frac{\prob_p({\cal J}_a^b)}{(G_{a-1}^{b-1})^2}
\geq C p \,e^{-C'p(b-a)^2},
$$
where $C,C'\in(0,\infty)$ depend only on $c_\pm$.
\end{lemma}

\begin{pf}
From the definition of ${\cal J}_a^b$ and Proposition \ref{recur} we
obtain
$$\prob_p({\cal J}_a^b) \geq [1-(1-p)^a]^4 (1-p)^{2b} \,p\, \exp \big[ -
(b-a)g(aq)-(b-a)g(bq)\big].$$ Note that $g$ is decreasing, and that
$(1-p)^k$ is bounded away from 0 and 1 for $k\in [c_-/p,c_+/p]$, so
we deduce
\begin{equation}\label{j-bound}
\prob_p({\cal J}_a^b) \geq C p\, \exp \big[-2(b-a)g(aq)\big].
\end{equation}
 Also we have
\begin{equation}\label{prod-bound}
(G_{a-1}^{b-1})^2=\exp \big[\textstyle - 2\sum_{i=a-1}^{b-2}
g(iq)\big] \leq \exp \big[-2(b-a)g(bq)\big].
\end{equation}
 Now $g(aq)-g(bq)\leq (bq-aq)
\max_{z\in[aq,bq]} |g'(z)|$, but the ratio $q/p$ is bounded for
$p<1/2$, hence $g'$ is uniformly bounded over the relevant interval,
and we obtain $g(aq)-g(bq)\leq C'(b-a)p$.  Therefore dividing
\eqref{j-bound} by \eqref{prod-bound} gives the result.
\end{pf}

\begin{pfof}{Proposition \ref{nuc-slow}}
Let $m=\lfloor Mp^{-1/2}\rfloor$, where $M<1/4$ is a constant to be
chosen later.  Suppose integers $(a_i,b_i)_{i=1,\ldots,m}$ and $B$
satisfy:
\begin{equation}\label{possibilities}
\begin{array}{c}
p^{-1}<a_1\leq b_1 \leq a_2 \leq \cdots \leq b_m<2 p^{-1}\leq B\\
b_i-a_i \in [4,p^{-1/2}] \qquad \forall i
\end{array}
\end{equation}

Let $C,C'$ be the constants from Lemma \ref{jig-bound} corresponding
to $c_-=1$ and $c_+=2$. Then from the definition of the event ${\cal
E}$ together with Lemmas \ref{properties}(i), \ref{diag-bound} and
\ref{jig-bound} we obtain:
\begin{align}
\prob_p\big[{\cal E}(a_1,\ldots,b_m)\big] & \geq p^4 \big[ C
p\,e^{-C'p(p^{-1/2})^2} \big]^m \exp \big[\textstyle -
2\sum_{i=1}^{B-1}
g(iq)\big] \nonumber \\
&= p^4 ( C'' p )^m (G_1^B)^2 \label{one-poss}
\end{align}
for $C''$ a fixed constant. Now since $m p^{-1/2}<p^{-1}/4$, the
number of possible choices of $(a_i,b_i)_{i=1,\ldots,m}$
satisfying \eqref{possibilities} is at least
\begin{equation}\label{num-poss}
{{\lfloor p^{-1} - mp^{-1/2} \rfloor} \choose m} (p^{-1/2}-4)^m
\geq \frac{(p^{-1}/2)^m}{m^m} (p^{-1/2}/2)^m = \Big(
\frac{1}{4pM}\Big)^m
\end{equation}
for $p$ sufficiently small.

By Lemma \ref{properties}(ii),(iii) we may multiply \eqref{one-poss}
and \eqref{num-poss} to give for $p$ sufficiently small and all
$B>2p^{-1}$,
$$I(B)\geq p^4 \Big( \frac{C''}{4M}\Big)^m (G_1^B)^2.$$
 Now choose $M=C''/8$ (recall that $C''$ was an absolute
constant) so that $C''/4M = 2$. Also note that since $g$ is
decreasing,
$$-\log G_1^B =\sum_{i=1}^{B-1} g(iq) \leq q^{-1} \int_0^{Bq} g \leq
p^{-1} \int_0^\infty g = p^{-1}\lambda.$$ Hence for $p$
sufficiently small,
$$I(B) \geq p^4 2^{Mp^{-1/2}/2} \exp[-2 p^{-1} \lambda]
\geq\exp[-2p^{-1}\lambda+cp^{-1/2}],$$
 as required.
\end{pfof}

\section{Explicit bound for the modified model}
\label{sec-mod}

In this section we prove Theorem \ref{explicit}.  Since we always
refer to the modified model we sometimes omit the subscript $M$ in
$I_M$.

\begin{prop}[nucleation centres]\label{mod-nuc}
Consider the modified model.  For any $p\leq 1/10$ and any $B\geq
\sqrt{2/p}$ we have
$$I(B) \geq \exp \Big[-2\lambda_M/q + 2\sqrt{2/p} -\log p^{-1} -3.2\Big],$$
where $\lambda_M=\pi^2/6$.
\end{prop}

\begin{lemma}[diagonal spanning]
\label{diag}  For the modified model we have for any positive
integer $a$ and any $p\in(0,1)$,
$$I_M(a)\geq\tfrac 12 \left(2p-p^2\right)^a.$$
\end{lemma}

\begin{pf}
Note that for $a\geq 2$, the square $R(a)$ is internally spanned
provided $(1,1)$ is occupied and $R(2,2;a,a)$ is internally
spanned, or alternatively provided $(1,a)$ is occupied and $R(2,1;
a,a-1)$ is internally spanned.  Hence
$$I(a)\geq pI(a-1)+(1-p)pI(a-1)= (2p-p^2)I(a-1).$$
The result follows by induction.
\end{pf}

Denote
$$F_a^b=F_a^b(p):=\textstyle \prod_{j=a}^{b-1} \big(1-(1-p)^j)=
\exp \Big[\textstyle  - \sum_{i=a}^{b-1} f(iq)\Big].$$
\begin{lemma}[growth]
\label{growth}
 Let $a\leq b$ be integers and let $p\in(0,1)$.  For
the standard or modified model, we have
$$I(b)\geq I(a) (F_a^b)^2.$$
\end{lemma}

\begin{pf}
Let $F$ be the event that each of the strips
\begin{align*}
R(j+1,1;\; j+1,j),\qquad j=a,a+1,\ldots,b, \\
R(1,j+1;\; j,j+1),\qquad j=a,a+1,\ldots,b
\end{align*}
is non-vacant. It is easily seen that if $R(a)$ is internally
spanned and $F$ occurs then $R(b)$ is internally spanned. Hence
$$I(b)\geq \prob(\{R(a) \text{ is i.s.}\}\cap F)=I(a)\prob(F)=I(a)(F_a^b)^2.$$
\end{pf}

We next note some elementary bounds.  We have
\begin{equation}
p\leq q\leq p+p^2, \label{bounds-q}
\end{equation}
where the second inequality holds provided $p<1/2$.
The function $F_a^b$ satisfies
\begin{equation}
\exp \bigg[-\frac 1q \int_{(a-1)q}^{(b-1)q} f \bigg]
\leq F_a^b \leq
\exp \bigg[-\frac 1q \int_{aq}^{bq} f\bigg], \label{trap}
\end{equation}
since $f$ is decreasing.

Also note the inequalitites
\begin{align}
\log \epsilon^{-1} &\leq f(\epsilon) \leq \log \epsilon^{-1} +\epsilon \label{f-small}
\\
e^{-K} &\leq f(K) \leq e^{-K}+e^{-2K}, \label{f-large}
\end{align}
where the fourth inequality holds provided $K>1/2$.
 (The inequalities are useful when $\epsilon\ll 1\ll K$).
Hence
\begin{align}
\epsilon \log \epsilon^{-1} +\epsilon  &\leq \int_0^\epsilon f \leq
\epsilon \log \epsilon^{-1} +\epsilon + \tfrac 12 \epsilon^2 \label{int-small}\\
e^{-K} &\leq \int_K^\infty f \leq e^{-K}+\tfrac 12 e^{-2K}, \label{int-large}
\end{align}
where the fourth inequality holds provided $K>1/2$.

\begin{pfof}{Proposition \ref{mod-nuc}}
Fix $p<1/10$, and let $A\le B$ be positive integers (later we
will take $A\approx \sqrt{2/p}$).

By Lemmas \ref{diag} and \ref{growth} we have
$$I(B)\geq \tfrac 12 (2p-p^2)^A (F_A^B)^2,$$
so using (\ref{trap}), (\ref{integral}) and (\ref{int-small}), and
rearranging,
\begin{align*}
\lefteqn{\log I(B)
  \geq -\log 2 + A \log(2p-p^2) - \frac{2}{q} \int_{(A-1)q}^\infty f} \\
 & \geq -\log 2 + A \log(2p-p^2) - \frac{2}{q}
 \Big( \lambda_M
 - (A-1)q \log [(A-1)q]^{-1} - (A-1)q \Big) \\
 & = -\frac{2 \lambda_M}{q}
 +2(A-1)\log \frac{e \sqrt{2}}{(A\!-\!1)\sqrt p}
+2(A-1)\log \frac{p}{q}
 + A\log (1-\frac p2) +\log p,
 \end{align*}
where we have written $(2p-p^2)=2p(1-p/2)$.
By (\ref{bounds-q}), for $p<1/2$ we have $\log (p/q)\geq \log
[p/(p+p^2)]=-\log(1+p)\geq -p$, and $\log(1-p/2)\geq -p/2-p^2/4$,
so we obtain
$$\log I(B)\geq
 -\frac{2 \lambda_M}{q}
 +2(A-1)\log \frac{e \sqrt{2}}{(A-1)\sqrt p}
-2(A-1)p
 - A(p/2+p^2/4) +\log p.$$

 Now let
$$A=\big\lceil \sqrt{2/p} \big\rceil,$$
to give
for $p\leq 1/10$ and $B\ge A$,
\begin{align*}
\lefteqn{\log I(B)} \\
&\geq
 -\frac{2 \lambda_M}{q}
 +2\big(\sqrt{2/p}-1\big)1  -2\sqrt{2/p} \;p
 - \big(\sqrt{2/p}+1)(p/2+p^2/4) +\log p  \\
 &\geq -\frac{2 \lambda_M}{q} +\frac{2\surd 2}{\surd p} -\log p^{-1} -
 3.2.
\end{align*}
Note the non-trivial cancelation between terms in $p^{-1/2}\log
p^{-1}$.
\end{pfof}

The following variant of Lemma \ref{comp}(i) allows better control
of the error terms.
\begin{lemma}[scanning estimate]
\label{scan}  Let $b,\ell,m$ positive integers with $mb<\ell$, and
let $p\in(0,1)$.  For the standard or modified model, we have
$$I(\ell)\geq\left(1-e^{-m^2I(b)}\right) (F_b^\ell
F_{\ell-mb}^\ell)^2 \big(1-(1-p)^{\ell-mb}\big)^\ell.$$
\end{lemma}

\begin{pf}
Consider the $m^2$ disjoint squares
$$S_k:=R(b)+bk,\qquad k\in\{0,\ldots,m-1\}^2,$$
and let
$$\{0,\ldots,m-1\big\}^2=\big\{k(1),k(2),\ldots,k(m^2)\big\}$$
be the lexicographic ordering of the set on the left side.
For $i=1,\ldots, m^2$ define the event
$$J_i=\{S_{k(i)} \text{ is internally spanned}\},$$
and let $F_i$ be the event that each of the strips
\begin{align*}
R(\ell)\cap [bk(i)+R(j+1,1;\; j+1,j)],\qquad j=b,b+1,\ldots \\
R(\ell)\cap [bk(i)+R(1,j+1;\; j,j+1)],\qquad j=b,b+1,\ldots
\end{align*}
that is non-empty is non-vacant. See Figure \ref{fig-scan}. Also
define the event
$$E=\big\{\cl{W\cap R(\ell)}\supseteq R(mb+1,mb+1;\;\ell,\ell)\big\}.$$
It is straightforward to see that for any $i$, if $J_i$ and $F_i$ occur then
$E$ occurs.  Furthermore, for each $i$,
the event $F_i$ is independent of the events $J_1,\ldots, J_i$.
\begin{figure}
\centering \resizebox{!}{2in}{\includegraphics{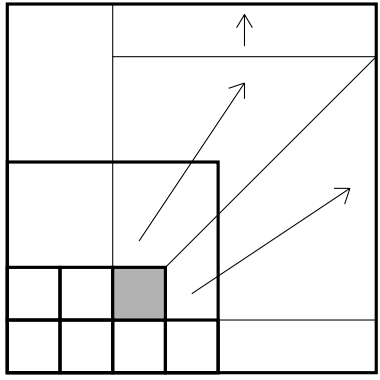}} \caption{An
illustration of the proof of Lemma \ref{scan}.  Here $m=4$, and
the first internally spanned sub-square is $S_{k(7)}=S_{(2,1)}$.
 The arrows indicate the event $F_7$.} \label{fig-scan}
\end{figure}

Hence we have
\begin{align}
\prob(E)& \geq \prob \bigg[ \bigcup_{i=1}^{m^2} \big(J_1^C \cap
\cdots\cap J_{i-1}^C\cap J_i
\cap F_i\big)\bigg] \nonumber\\
&=  \sum_{i=1}^{m^2} \prob \big(J_1^C \cap \cdots\cap J_{i-1}^C\cap J_i\big)
\prob(F_i) \nonumber\\
&\geq  \prob(J_1 \cup\cdots \cup J_{m^2}) \min_i \prob(F_i) \nonumber\\
&\geq  \label{scan-step} \left(1-e^{-m^2I(b)}\right)(F_b^\ell)^2
\big(1-(1-p)^{\ell-mb}\big)^\ell.
\end{align}

To conclude, let $H$ be the event that each of the strips
\begin{align*}
R(j,j-1;\; j,\ell),\qquad j=mb,\ldots,2,1 \\
R(j-1,j;\; \ell,j),\qquad j=mb,\ldots,2,1
\end{align*}
is non-vacant.  Using the Harris-FKG inequality we have $I(\ell)\geq
\prob(E\cap H)\geq \prob(E)\prob(H) \geq
\prob(E)(F_{\ell-mb}^\ell)^2$, and combining this with
(\ref{scan-step}) gives the result.
\end{pf}

\begin{pfof}{Theorem \ref{explicit}}
Fix $p\leq 10$ and let $B\geq \sqrt{2/p}$, and take $L$ and $m$ such
that $L\geq mB$.  We use Lemma \ref{scan} to derive a lower bound
for $I(L)$.  We obtain
 \begin{equation}
I(L)\geq \Big( 1-e^{-m^2 I(B)} \Big) (F_B^\infty F_{L-mB}^\infty)^2
e^{-L f([L-mB]q)}. \label{bound-il}
\end{equation}
 Consider the first factor above.  Take
\begin{equation}
\label{def-m} m=\left\lceil \exp\bigg( \frac{\lambda_M}{q} -
{\frac{\surd 2}{\surd p}}
 + \frac{1}{2}\log p^{-1} + 1.8
 \bigg) \right\rceil.
\end{equation}
 Then Proposition \ref{mod-nuc} implies $\log (m^2 I(B)) \geq 0.4$, and therefore
$$1-e^{-m^2 I(B)} \geq 1-e^{-e^{0.4}}.$$

Turning to the other factors in (\ref{bound-il}), we have by
(\ref{trap}),
\begin{eqnarray*}
\lefteqn{(F_B^\infty F_{L-mB}^\infty)^2 e^{-L f([L-mB]q)}}\\
 &\geq& \exp\Big(-\frac{2}{q}\int_{(B-1)q}^\infty f
 - \frac{2}{q}\int_{(L-mB-1)q}^\infty f
 - L f([L-mb]q) \Big) \\
 &\geq& 1 -\frac{2}{q}\int_{(B-1)q}^\infty f
 - \frac{2}{q}\int_{(L-mB-1)q}^\infty f
 - L f([L-mb]q).
\end{eqnarray*}
We now set
\begin{equation}
\label{def-bl}
B=1+\Big\lceil\frac{3+\log q^{-1}}{q}\Big\rceil
 \quad\text{and}\quad
L=mB+4c q^{-2},
\end{equation}
for any $c\geq 1$.
It is straightforward to check
that for $p\leq 1/10$ we have $(L-mB-1)q>(B-1)q>1/2$, so we may use
(\ref{f-large}),(\ref{int-large}) to bound the above terms thus:
$$\frac{2}{q}\int_{(B-1)q}^\infty f
 - \frac{2}{q}\int_{(L-mB-1)q}^\infty f
 \leq \frac{4}{q}\Big(e^{-(B-1)q}+e^{-2(B-1)q}\Big) \leq 4 e^{-3}+4 e^{-6},$$
 and
\begin{multline*}
 L f([L-mb]q)\leq 2L e^{-(L-mB)q}\leq
2\big(e^{2/q}2q^{-2}+4cq^{-2}+1)e^{-4c/q} \\
\leq 2(e^{2/q}2q^{-2}+4q^{-2}+1) e^{-4/q}\leq e^{-2}
\end{multline*}
since $m\leq e^{2/q}$ and $B\leq 2q^{-2}$ for $p\leq 1/10$. Hence,
returning to (\ref{bound-il}), for the given choices of $B,L$ we
have
$$I(L)\geq (1-e^{-e^{0.4}})(1-4e^{-3}-4e^{-6}-e^{-2})>1/2.$$

From (\ref{def-bl})
we have shown that $I(L,p)>1/2$ provided $p\leq1/10$ and
\begin{equation}
\label{bound-plogl}
p\log L \geq  p\log(mB+4q^{-2}) =p\log m+p\log B +p\log\Big( 1+\frac{4q^{-2}}{mB}\Big).
\end{equation}
Finally we need to find upper bounds for the terms appearing on the
right of \eqref{bound-plogl}. By \eqref{def-m} we have
$$p\log m \leq \lambda_M\frac{p}{q} -\sqrt{2p} +\frac{1}{2} p\log p^{-1}
+1.8p+p\log\frac{m}{m-1}.$$
But for $p\leq 1/10$ we have
$p\log(m/(m-1))=-p\log(1-1/m)\leq 2p/m\leq 2pe^{-1/p}\leq 0.001 p$, while
$p/q\leq p/(p+p^2/2)\leq 1-0.47 p$, so
$$p\log m \leq \lambda_M-\sqrt{2p}+\frac{1}{2} p\log p^{-1}
+1.03p.$$
By \eqref{def-bl} we have
$$p\log B \leq p\log\Big(2+\frac{3}{q}+\frac{\log q^{-1}}{q}\Big)\leq p\log(2.6 p^{-1.3})
=0.96 p+1.3p\log p^{-1}.$$
Since $4q^{-2}>B$ and $m\geq e^{1/p}$ for $p\leq 1/10$, we have
$$p\log\Big( 1+\frac{4q^{-2}}{mB}\Big)\leq p e^{-1/p}\leq 0.001 p.$$
Hence the right side of \eqref{bound-plogl} is at most
$$\lambda_M-\sqrt{2p}+1.8 p\log p^{-1}
+ 2 p,$$ as required.
\end{pfof}

\section*{Open Problems}

\begin{mylist}
\item Prove a complementary bound to Theorem \ref{slow}.  For example, do there exist
$\gamma,c<(0,\infty)$ such that $(L,p)\to(\infty,0)$ with $p\log
L<\lambda -c(\log L)^{-\gamma}$ implies $I\to 0$?
\item Prove matching upper and lower bounds, e.g. involving inequalities of the form
$p\log L \lessgtr \lambda-c(\log L)^{\gamma\pm\epsilon}$, or even
$p\log L \lessgtr \lambda-(c\pm \epsilon)F(L)$.
\item Extend the results to other bootstrap percolation models
for which sharp thresholds are known to exist -- currently those in \cite{h-mod,h-l-r}.
\item Identify more precisely the width of the critical window as $p$
varies.  Is it the case that $p_{1-\epsilon}\log L-p_\epsilon\log
L=\Theta(1/\log L)$ as $L\to\infty$?
\end{mylist}

\bibliography{bib}

\section*{}
\vspace{-10mm} \sc

\noindent Alexander E. Holroyd: {\tt holroyd(at)math.ubc.ca}\\
Department of Mathematics, University of British Columbia,\\
121-1984 Mathematics Rd, Vancouver, BC V6T 1Z2, Canada \\[2mm]

\noindent Janko Gravner: {\tt gravner(at)math.ucdavis.edu}\\
Mathematics Department, University of California,\\
Davis, CA 95616, USA

\end{document}